\documentclass[10pt,reqno]{amsart}
\usepackage{geometry} 
\usepackage{graphicx}
\usepackage{amsmath}
\usepackage{color}

\newtheorem{thm}[equation]{Theorem}
\newtheorem{pro}[equation]{Proposition}
\newtheorem{cor}[equation]{Corollary}
\newtheorem{lem}[equation]{Lemma}

\theoremstyle{definition}

\newtheorem{DEF}[equation]{Definition}
\newtheorem{rem}[equation]{Remark}

\usepackage{amsfonts,amssymb,amscd,amsmath,enumerate,verbatim,calc}

\def\span{\hbox{span}}

\def\bbbf{\mathbb{F}}

\def\ss{{\mathcal S}}

\def\lam{\lambda}
\def\Lam{\Lambda}

\def\qed{\hfill$\Box$\vspace{3mm}}

\def\u{{\mathcal U}}

\def\pp{{\mathcal P}}
\def\vv{{\mathcal V}}
\def\uu{{\mathcal U}}

\def\lam{\Lambda_I^{\bar{i}}}

\def\lam{\lambda}
\def\Lam{\Lambda}
\def\uu{{\mathcal U}}
\def\p{{ P}}
\begin{document}

\markboth{Non-commutative Poisson  algebras with a set grading} {V. Khalili}

\date{}

\centerline{\bf Non-commutative Poisson  algebras with a set grading }

\vspace{.5cm}\centerline{Valiollah Khalili\footnote[1]{Department
of mathematics, Faculty of sciences, Arak University, Arak 385156-8-8349, Po.Box: 879, Iran. 
 V-Khalili@araku.ac.ir\\
\hphantom{ppp}2000 Mathematics Subject Classification(s): 17A30,
17B63, 17A60.
\\
\hphantom{ppp}Keywords: Non-commutative Poisson  algebras, set-graded Lie  algebra, structure theory }\;\;}

\vspace{1cm} \noindent ABSTRACT: In this paper we   study of the structure of  non-commutative Poisson  algebras with an arbitrary set $\ss.$  
We show that any of such an algebra $\pp$ decomposes as
$$\pp=\uu\oplus\sum_{[\lambda]\in(\Lambda_\ss\setminus\{0\})/\sim}\pp_{[\lambda]},$$
where $\uu$ is a linear subspace complement of
$\span_{\bbbf}\{ [\pp_{\mu}, \pp_{\eta}]+\pp_{\mu}\pp_{\eta}  : \mu, \eta\in[\lam]\}\cap\pp_0$ in $\pp_0$ and any $\pp_{[\lambda]}$ a well-described
graded ideals of $\pp,$ satisfying $[\pp_{[\lambda]},
\pp_{[\mu]}]+\pp_{[\lambda]}
\pp_{[\mu]}=0$ if $[\lambda]\neq[\mu].$
 Under certain conditions, the simplicity of $\pp$ is
characterized  and it is shown that $\pp$  is the direct sum  of the family of its graded simple ideals.

\vspace{1cm} \setcounter{section}{0}
\section{Introduction}\label{introduction}

A non-commutative Poisson algebra is a Lie algebra endowed with a non-commutative associative product in such a way that the Lie and associative products are compatible via a Leibniz rule, which initially  appeared in  the work of Smeon-Denis Poisson in the theory of celestial mechanics \cite{A}. These algebras play a central role in the
study of Poisson geometry \cite{A1, L, V}, in deformation quantization \cite{K, H} and  in deformation of commutative associative algebras \cite{G}. Poisson  structures are also used in the
study of vertex operator algebras (see \cite{FB}). Moreover, the cohomology group, deformation, tensor product
and $\Gamma$-graded of Poisson algebras have been studied by many authors in \cite{AD, GK, SS, WGZ, WZC, ZL}.

 There is a concept of grading not involving groups. This is a grading by means of an arbitrary set $\ss,$ not necessarily a group. The study of gradings on Lie algebras \cite{BSZ, E, EK},  begins in the 1933 seminal Jordan’s work,
with the purpose of formalizing Quantum Mechanics \cite{J}. Set gradings on Lie algebras where
first considered in the literature by Patera and Zassenhaus in \cite{PZ} under the restriction that the set $\ss$ agrees with the support of the grading which they called these gradings Lie gradings.  It is worth mentioning that the so-called techniques of connection of roots had long
been introduced by Calderon, Antonio J, on split Lie algebras with symmetric root
systems in \cite{C1}.  Since then, the interest on gradings by the technique of connections of elements in the support of the
graing on different classes of algebras has been remarkable in the recent years, motivated in part by their application in physics and geometry. Recently, in \cite{ABCS, C2, Kh1},  Lie algebras with a set grading,  Lie superalgebras with a set grading and  Leibniz superalgebras with a set graing,    have been determined
by the connections of the support of the grading. Our goal in this work is to study the structure of arbitrary non-commutative Poisson  algebras (not necessarily
 finite-dimensional)  over an arbitrary base field $\bbbf$ by focusing on the technique of connections of elements in the support of the grading.  their structure. Inspired by the previous works on split non-commutative Poison algebras in \cite{C0}, we would like to 
 study the set-graded non-commutative Poisson algebras by focusing on their inner structure by using of the  connection techniques.

Throughout this paper, non-commutative Poisson  algebras $\pp$ are considered of arbitrary dimension and over an
arbitrary base field  $\bbbf,$ with characteristic zero.  We also
consider an arbitrary set  $\ss$ with identity zero.

To close this introduction, we briefly outline the contents of
the paper. In Section 2, we begin by recalling the necessary
background  on graded Poisson color algebras.
Section 3  develops the technique of connections of elements in the support of the
grading for graded non-commutative Poisson  algebras by means of an arbitrary set. We also show
that such an  algebra $\pp$ with a   support $\Lam_{\ss}$  is of the form
$\pp=\uu\oplus\sum_{[\lambda]\in(\Lambda_\ss\setminus\{0\})/\sim}\pp_{[\lambda]},$ with $\uu$   a vector subspace complement  of $\span_{\bbbf}\{ [\pp_{\mu}, \pp_{\eta}]+\pp_{\mu}\pp_{\eta}  : \mu, \eta\in[\lam]\}\cap\pp_0$ in $\pp_0$ and any $\pp_{[\lambda]}$ a well-described
graded ideals of $\pp,$ satisfying $[\pp_{[\lambda]},
\pp_{[\mu]}]+\pp_{[\lambda]}
\pp_{[\mu]}=0$ if $[\lambda]\neq[\mu].$  In
section 4, we show that  under certain conditions, in the case of
$\pp$ being of maximal length, the gr-simplicity of the algebra is characterized.\\

\section{Preliminaries} \setcounter{equation}{0}\

In this section, we first recall some basic definitions and properties of set-graded non-commutative Poisson algebras.

\begin{DEF}\cite{C0} A {\em   non-commutative Poisson algebra}   is a Lie algebra $(\pp, [., .])$ over
an arbitrary base field $\bbbf,$ endowed with an associative product, denoted by
juxtaposition, in such a way the following Leibniz  rule:
$$
[x, y z]=[x, y]z+y[x, z],
$$
holds for any  $x, y, z\in\pp.$
\end{DEF}
A subalgebra $A$   is a linear subspace of $\pp$ closed by both the Lie and the
associative products, that is $[A, A]+AA\subset A$. An ideal $I$ of $\pp$  is a subalgebra satisfying $[I, \pp]+I\pp+\pp I\subset I.$
A non-commutative Poisson algebra $\pp$  is called simple if $[\pp, \pp]\neq 0, \pp\pp\neq 0$
 and its only ideals are $\{0\}$ and $\pp.$

\begin{DEF}\label{set-grPoisson} Let
$\pp$ be a  non-commutative Poisson algebra and
$\ss$ be an arbitrary (non-empty) set. It is said that $\pp$ has
a {\em set grading by means of $\ss$} or it is {\em set-graded} if
$$
\pp=\bigoplus_{s\in\ss}\pp_s,
$$
as a direct sum of non-zero subspaces indexed by $\ss,$   having the
property that, for any $s, t\in\ss$ with $0\neq[\pp_s, \pp_t]+\pp_s \pp_t,$
there exists unique $u\in\ss$ such that $[\pp_s, \pp_t]+\pp_s \pp_t\subset\pp_u.$ 
\end{DEF}

Note that, split non-commutative Poison algebras, graded Poisson color algebras are examples of set-graded Poisson superalgebras. Hence, the present paper  extends the results in \cite{C0, Kh2}.

We call the {\em support} of the set-grading to be the set
$$
\Lambda_\ss=\{\lambda\in\ss: \pp_\lambda\neq0\}.
$$
So we can write
$$
\pp=\bigoplus_{\lambda\in\Lambda_\ss}\pp_\lambda.
$$

The usual regularity conditions will be understood in the graded
sense compatible with the non-commutative  Poison algebra structure, that is, a
{\em graded subalgebra} $A$ of a set-graded  non-commutative  Poison algebra $\pp$
is a graded subspace such that it  splits as
\begin{equation}\label{0}
A=\bigoplus_{s\in\ss}A_s~~~\hbox{with}~~A_s=A\cap\pp_s,
\end{equation}
and that $[A, A]+A A\subset A.$ A  graded subalgebra $I$ of $\pp$ is a
{\em graded ideal} if $[I, \pp]+I\pp+\pp I\subset I.$    A set-graded
 non-commutative  Poison algebra $\pp$ will be called {\em gr-simple} if $[\pp,
\pp]+\pp \pp\neq0$ and its only ideals are $\{0\}$ and $\pp.$\\

\section{Connections in $\Lambda_\ss. Decompositions$} \setcounter{equation}{0}\

In this section, we begin by developing connection techniques  in the support of a set-graded  non-commutative  Poisson algebras as in \cite{C2}.

Let $\pp=\bigoplus_{\lam\in\Lam_\ss}\pp_\lambda$ be an arbitrary  set-graded non-commutativ Poisson algebra,  with $\Lam_\ss$
the support of the grading. First we recall some terminology
which can be found in \cite{C2}.

For each $\lambda\in\Lambda_\ss,$ a new symbol
$\tilde{\lambda}\notin\Lambda_\ss$ is introduced and we denote by
$$
\widetilde{\Lambda_\ss} : =\{\tilde{\lambda} :
\lambda\in\Lambda_\ss\},
$$
the set consisting of all these new symbols.

We will denote by $\p(A)$ the power set of a given set $A$. Next, we consider the following operation
$$
\star :
(\Lambda_\ss\dot{\cup}\widetilde{\Lambda_\ss})\times(\Lambda_\ss\dot{\cup}\widetilde{\Lambda_\ss})\longrightarrow\p(\Lambda_\ss),
$$
given by

$\bullet$ for $\lambda, \mu\in\Lambda_\ss$,
$$
\lambda\star\mu=\left\{\begin{array}{ll}
\emptyset&\hbox{~~if~~}[\pp_\lam, \pp_\mu]+\pp_\lam \pp_\mu=0\\
\{\eta\}&\hbox{~~if~~} 0\neq[\pp_\lam, \pp_\mu]+\pp_\lam \pp_\mu\subset\pp_\eta
\end{array}\right.~;
$$

$\bullet$  for $\lambda\in\Lambda_\ss,
\tilde{\mu}\in\tilde{\Lambda_\ss}$,
$$
\lambda\star\tilde{\mu}=\tilde{\mu}\star\lambda=\{\eta\in\Lambda_\ss
: 0\neq [\pp_\eta, \pp_\mu]+\pp_\eta \pp_\mu\subset\pp_\lambda \};
$$

$\bullet$  for $\tilde{\lambda},
\tilde{\mu}\in\tilde{\Lambda_\ss}$,
$$
\tilde{\lambda}\star\tilde{\mu}=\emptyset.
$$
From now on, given any $\tilde{\lambda}\in\Lambda_\ss$ we will
denote
$$
\tilde{\tilde{\lambda}} :=\lambda.
$$
Given also any subset $\Omega$ of
$\Lambda_\ss\dot{\cup}\widetilde{\Lambda_\ss},$ we write
$\tilde{\Omega} :=\{\tilde{\lambda} : \lambda\in\Omega\}$ if
$\Omega\neq\emptyset$ and also $\tilde{\emptyset} :=\emptyset.$

It is worth to mention that,  sometimes it is interesting to distinguish one element $0$ in the support of the grading,  because the   homogeneous space $\pp_0$  has, in a sense,
a special behavior to the remaining elements in the set of spaces $\pp_\lam,$ for $\lam\in\Lam_\ss.$ This is for instance 
case in which the grading set $\ss$ is an abelian group, where
the homogeneous space $\pp_0$ associated to the unit element $0$
in the group enjoys a distinguished role (see \cite{Kh2}).  If we consider
the group-grading determined by the Cartan decomposition of a
split Lie algebra of maximal length, the  homogeneous
space associated to the unit element agrees with the Cartan
subalgebra $H$ and being $\dim\pp_\lambda=1$ for
any $\lambda$ in the support of the grading up to $\dim\pp_0$
which is not bounded by this condition (see \cite{C1}). From here, we are going
to feel free in our study to distinguish one element $0$ in the
support of the grading. Hence, let us now fix an element $0$ such
that either $0\in\Lambda_\ss$ satisfying the property
$0\star\lambda\neq\{0\}$ for any
$\lambda\in\Lambda_\ss\setminus\{0\},$ or $0=\emptyset.$ The
possibility $0=\emptyset$ holds for the case in which it is not wished
to distinguish  any element in $\Lambda_\ss.$ Finally, we need to introduce the following mapping;
$$
\psi :
\left((\Lambda_\ss\dot{\cup}\widetilde{\Lambda_\ss})\setminus\{0,
\tilde{0}
\}\right)\times(\Lambda_\ss\dot{\cup}\widetilde{\Lambda_\ss})\longrightarrow\p\left((\Lambda_\ss\dot{\cup}\widetilde{\Lambda_\ss})\setminus\{0,
\tilde{0} \}\right),
$$
given by

$\bullet~~ \psi(\emptyset,
\Lambda_\ss\dot{\cup}\widetilde{\Lambda_\ss})=\emptyset$;

$\bullet$  for  any
$\emptyset\neq\Omega\in\p\left((\Lambda_\ss\dot{\cup}\widetilde{\Lambda_\ss})\setminus\{0,
\tilde{0} \}\right)$ and any
$a\in(\Lambda_\ss\dot{\cup}\widetilde{\Lambda_\ss})$,
$$
\psi(\Omega,a)=\left(\left (\bigcup_{x\in\Omega}(x\star
a)\right)\setminus\{0\} \right) \cup\left(\left(\bigcup_{x\in\Omega}\widetilde{(x\star
a)}\right)\setminus\{0\}\right).
$$
\begin{rem}\label{1111}
It is obvious that for any
$\Omega\in\p\left((\Lambda_\ss\dot{\cup}\widetilde{\Lambda_\ss})\setminus\{0,
\tilde{0} \}\right)$ and any
$a\in(\Lambda_\ss\dot{\cup}\widetilde{\Lambda_\ss})$, we have
\begin{equation}\label{tild}
\psi(\Omega, a)=\widetilde{\psi(\Omega, a)},
\end{equation}
and
$$
\psi(\Omega, a)\cap\Lambda_\ss=\left(\bigcup_{x\in\Omega}(x\star
a)\right)\setminus\{0\}.
$$
Also observe that for any $\lambda\in\Lambda_\ss$ and
$a\in(\Lambda_\ss\dot{\cup}\widetilde{\Lambda_\ss})$ we have
$\lambda\in x\star a$ for some $x\in \Lambda_\ss$ if and only if
$x\in\lambda\star\tilde{a},$ while $\lambda\in\mu\star a$ for
some $\mu\in\widetilde{\Lambda_\ss}$ if and only if
$\tilde{\mu}\in\tilde{\lambda}\star a.$ These facts together with
Eq.(\ref{tild}) imply that for any
$\Omega\in\p\left((\Lambda_\ss\dot{\cup}\widetilde{\Lambda_\ss})\setminus\{0,
\tilde{0} \}\right)$ such that $\Omega=\tilde{\Omega}$ and
$a\in(\Lambda_\ss\dot{\cup}\widetilde{\Lambda_\ss})$  we have
$\lambda\in \psi(\Omega, a)\cap\Lambda_\ss$ if and only if
$\lambda\in\Lambda_\ss$ and either $\psi(\{\lambda\},
\tilde{a})\cap\Omega\cap\Lambda_\ss\neq\emptyset$ or
$\psi(\widetilde{\Lambda_\ss},
a)\cap\Omega\cap\Lambda_\ss\neq\emptyset.$
\end{rem}

\begin{DEF}\label{conn}
Let $\lambda, ~\mu\in\Lambda_\ss\setminus\{0\}.$ We say that
$\lambda$ is connected to $\mu$ and denote it by $\lambda\sim\mu$
if there exists a family
$$
\{\lam_1, \lam_2, \lam_3, ...,
\lam_n\}\subset\Lambda_\ss\dot{\cup}\widetilde{\Lambda_\ss},
$$
satisfying the following conditions;\\

If $n=1:$
\begin{itemize}
\item[(1)] $\lam_1=\lambda=\mu.$
\end{itemize}

If $n\geq2:$
\begin{itemize}
\item[(1)] $\lam_1\in\{\lambda, \tilde{\lambda}\}$,
\item[(2)] $\psi(\{\lam_1\}, \lam_2)\neq\emptyset,\\
\psi(\psi(\{\lam_1\}, \lam_2), \lam_3)\neq\emptyset,\\
\psi(\psi(\psi(\{\lam_1\}, \lam_2), \lam_3), \lam_4)\neq\emptyset,\\
...\\
\psi(\psi(...(\psi(\{\lam_1\}, \lam_2), ...), \lam_{n-2}),
\lam_{n-1})\neq\emptyset$.\\
\item[(3)] $\mu\in\psi(\psi(...(\psi(\{\lam_1\}, \lam_2), ...), \lam_{n-2}),
\lam_{n-1}), \lam_n)$.
\end{itemize}
The family $\{\lam_1, \lam_2, \lam_3, ...,
\lam_n\}$ is
called a {\em connection} from $\lambda$ to $\mu.$
\end{DEF}

\begin{pro}\label{rel}
The relation $\sim$ in $\Lambda_\ss\setminus\{0\},$ defined by
$\lambda\sim\mu$ if and only if $\lambda$ is connected to $\mu,$
is an equivalence relation.
\end{pro}
\noindent {\bf Proof.} The proof is similar to the
proof  of Proposition 2.1 in \cite{C2}. \qed\\

By the above proposition, given $\lam\in\Lambda_\ss\setminus\{0\},$ we can consider the quotient set
$$
(\Lambda_\ss\setminus\{0\})/\sim :=\{[\lam]~:~\lam\in\Lambda_\ss\setminus\{0\}\},
$$
then $[\lam]$  is the set of elements in  in
$\Lambda_\ss\setminus\{0\},$ which are connected to $\lambda.$ By Proposition \ref{rel}, if $\mu\notin[\lam]$ then $[\lam]\cap[\mu]=\emptyset.$

Our next goal in this section is to associate an ideal
$I_{[\lambda]}$ of $\pp$ to any $[\lambda].$ Fix
$\lambda\in\Lambda_\ss\setminus\{0\},$ we  define the
set $I_{0, [\lambda]}$ by;
\begin{equation}\label{2.1}
\nonumber I_{0, [\lambda]}=\span_{\bbbf}\{ [\pp_{\mu}, \pp_{\eta}]+\pp_{\mu}\pp_{\eta}  : \mu, \eta\in[\lam]\}\cap\pp_0.
\end{equation}
That is
\begin{equation}\label{2.2}
I_{0, [\lam]}=\left(\sum_{\substack{\mu, \eta\in[\lam]\\ \mu\star\eta=\{0\}}}([\pp_{\mu}, \pp_{\eta}]+\pp_{\mu}\pp_{\eta})\right)\subset\pp_{ 0},
\end{equation}
where $\pp_0 :=\{0\}$ whence $0=\emptyset.$
Next, we define
$$
\vv_{[\lam]}
:=\bigoplus_{\mu\in[\lam]}\pp_\mu.
$$
Finally, we denote by $I_{[\lam]}$ the direct sum of the two graded
subspaces above, that is,
$$
I_{[\lam]} :=I_{0, [\lam]}\oplus\vv_{[\lam]}.
$$

\begin{pro}\label{subalg} For any
$\lambda\in\Lambda_\ss\setminus\{0\},$ the graded subspace
$I_{[\lambda]}$ is a graded subalgebra of $\pp.$
\end{pro}
\noindent {\bf Proof.}   First, we are going to check that
$I_{[\lam]}$ satisfies $[I_{[\lam]}, I_{[\lam]}]\subset I_{[\lam]}.$ We have
\begin{eqnarray}\label{200}
\nonumber[I_{[\lam]}, I_{[\lam]}]&=&[I_{ 0, [\lam]}\oplus\vv_{[\lam]},
I_{0,  [\lam]}\oplus\vv_{[\lam]}]\\
&\subset&[I_{ 0, [\lam]}, I_{ 0, [\lam]}]+[I_{ 0, [\lam]}, \vv_{[\lam]}]+[\vv_{[\lam]}, I_{ 0, [\lam]}]+[\vv_{[\lam]}, \vv_{[\lam]}].
\end{eqnarray}
Let us consider the second summand in (\ref{200}). Taking into account $I_{ 0, [\lam]}\subset\pp_0$  and suppose that there exists $\mu\in[\lam]$ such that  $0\neq[I_{ 0, [\lam]}, \pp_\mu]\subset[\pp_0, \pp_\mu]\subset\pp_\tau$ with  $0\star\mu=\{\tau\},~~\tau\in\Lambda_\ss\setminus\{0\}.$ Then the connection $\{\tau, 0\}$ gives us $\mu\sim\tau$ and so $\tau\in[\lam].$ Hence, $[I_{ 0, [\lam]}, \vv_{[\lam]}]\subset\vv_{[\lam]}.$ By  skew-symmetry of the Lie product $[\vv_{[\lam]}, I_{ 0, [\lam]}]\subset\vv_{[\lam]}.$ Therefore,

\begin{equation}\label{100}
[I_{ 0, [\lam]}, \vv_{[\lam]}]+[\vv_{[\lam]}, I_{ 0, [\lam]}]\subset I_{[\lambda]}.
\end{equation}

Consider now  the  fourth  summand in (\ref{200}). Suppose there exist $\mu, \eta\in[\lam]$ with $\mu\star\eta=\{\tau\},~\tau\in\Lam_{\ss}$ such that $0\neq[\pp_{\mu}, \pp_{\eta}]\subset\pp_\tau.$ If $\tau=0,$  clearly
$[\pp_\mu, \pp_\eta]\subset\pp_0,$ and taking into acount Eq. (\ref{2.2}) we have $[\pp_\mu, \pp_\eta]\subset I_{ 0, [\lam]}.$ Otherwise, if  $\tau\in\Lambda_\ss\setminus\{0\},$ then the connection $\{\tau, \tilde\eta\}$ gives us $\mu\sim\tau$ and so $\tau\in[\lam].$  Hence $[\pp_\mu, \pp_\eta]\subset\pp_{\tau}\subset\vv_{[\lam]}.$ In any case, we have
\begin{equation}\label{203}
[\vv_{[\lam]}, \vv_{[\lam]}]\subset I_{[\lam]}.
\end{equation}

Finally consider  the first summand in (\ref{200}), we have
\begin{eqnarray}\label{204}
\nonumber[I_{ 0, [\lam]}, I_{0,  [\lam]}]&=&\left[\sum_{\substack{\mu, \eta\in[\lam]\\ \mu\star\eta=\{0\}}}\left([\pp_{\mu}, \pp_{\eta}]+\pp_{\mu} \pp_{\eta}\right), \sum_{\substack{\mu', \eta'\in [\lam]\\ \mu'\star\eta'=\{0\}}}\left([\pp_{\mu'}, \pp_{\eta'}]+\pp_{\mu'} \pp_{\eta'}\right)\right]\\
\nonumber&\subset&\sum_{\substack{\mu, \eta, \mu', \eta'\in [\lam]\\ \mu\star\eta=\{0\}=\mu'\star\eta'}}(\left[ [\pp_{\mu}, \pp_{\eta}], [\pp_{\mu'}, \pp_{\eta'}]\right]+\left[ [\pp_{\mu}, \pp_{\eta}], \pp_{\mu'} \pp_{\eta'}\right]\\
&+&[ \pp_\mu \pp_\eta, [\pp_{\mu'}, \pp_{\eta'}]]+[\pp_\mu \pp_\eta, \pp_{\mu'} \pp_{\eta'}]).
\end{eqnarray}
For the first item in (\ref{204}), suppose there exist $\mu, \eta, \mu', \eta'\in [\lam]$ with $\mu\star\eta=\{0\}$ and $\mu'\star\eta'=\{0\}$ such that $[ [\pp_{\mu}, \pp_{\eta}], [\pp_{\mu'}, \pp_{\eta'}]]\neq 0.$ Taking into account Eq. (\ref{2.2}) and the Jacobi identity, we have
\begin{eqnarray*}
\left[ [\pp_{\mu}, \pp_{\eta}], [\pp_{\mu'}, \pp_{\eta'}]\right]&\subset&\left[\left[ [\pp_{\mu}, \pp_{\eta}], \pp_{\mu'}\right], \pp_{\eta'}\right]+\left[\left[ [\pp_{\mu}, \pp_{\eta}], \pp_{\eta'}\right], \pp_{\mu'}\right]\\
&\subset&\left[[\pp_0, \pp_{\mu'}],  \pp_{\eta'}\right]+\left[[\pp_0,  \pp_{\eta'}], \pp_{\mu'}\right].
\end{eqnarray*}
If $0\neq[\pp_0, \pp_{\mu'}]\subset\pp_\tau,$ with $0\star\mu'=\{\tau\},~\tau\in\Lambda_\ss\setminus\{0\}$ then $\{\tau, 0\}$ is a connection from $\tau$ to $\mu',$ so $\tau\in[\lam].$ Similarly, with $[\pp_0,  \pp_{\eta'}]$ to get

\begin{equation}\label{204.1}
\sum_{\substack{\mu, \eta, \mu', \eta'\in [\lam]\\ \mu\star\eta=\{0\}=\mu'\star\eta'}}[ [\pp_{\mu}, \pp_{\eta}], [\pp_{\mu'}, \pp_{\eta'}]]\subset I_{0,  [\lam]}.
\end{equation}

For the second and  third items in (\ref{204}), taking into account $\pp_{\mu'}\pp_{\eta'}\subset\pp_0$  or $\pp_{\mu}\pp_{\eta}\subset\pp_0$ and  the  Jacobi identity, one can get
\begin{equation}\label{204.2}
\sum_{\substack{\mu, \eta, \mu', \eta'\in [\lam]\\ \mu\star\eta=\{0\}=\mu'\star\eta'}}([ [\pp_{\mu}, \pp_{\eta}], \pp_{\mu'} \pp_{\eta'}]+[ \pp_{\mu} \pp_{\eta}, [\pp_{\mu'}, \pp_{\eta'}]])\subset I_{0,  [\lam]}.
\end{equation}

For the last item in (\ref{204}), suppose there exist $\mu, \eta, \mu', \eta'\in [\lam]$ with $\mu\star\eta=\{0\}$ and $\mu'\star\eta'=\{0\}$ such that $[\pp_{\mu} \pp_{\eta}, \pp_{\mu'} \pp_{\eta'}]\neq 0.$ Taking into account Eq. (\ref{2.2}) and   the Leibniz rule, we have

\begin{eqnarray*}\label{204.3}
[\pp_{\mu} \pp_{\eta}, \pp_{\mu'} \pp_{\eta'}]&\subset&[\pp_0, \pp_{\mu'} \pp_{\eta'}]\\
&\subset&[\pp_0, \pp_{\mu'}] \pp_{\eta'}+\pp_{\mu'}[\pp_0,  \pp_{\eta'}])
\end{eqnarray*}
Again as above, If $0\neq[\pp_0, \pp_{\mu'}]\subset\pp_\tau,$ with $0\star\mu'=\{\tau\},~\tau\in\Lambda_\ss\setminus\{0\}$ then $\{\tau, 0\}$ is a connection from $\tau$ to $\mu',$ so $\tau\in[\lam].$ Similarly, with $[\pp_0,  \pp_{\eta'}]$ to get
\begin{equation}\label{204.4}
\sum_{\substack{\mu, \eta, \mu', \eta'\in [\lam]\\ \mu\star\eta=\{0\}=\mu'\star\eta'}}[\pp_{\mu} \pp_{\eta}, \pp_{\mu'} \pp_{\eta'}]\subset I_{0,  [\lam]}.
\end{equation}

From Eqs. (\ref{204.1})-(\ref{204.4}),   we conclude that

\begin{equation}\label{205}
[I_{0, [\lam]}, I_{0, [\lam]}]\subset I_{0, [\lam]}.
\end{equation} 
Finally,  Eqs (\ref{100}), (\ref{203}) and (\ref{205}) give us 
\begin{equation}\label{2061}
[I_{[\lam]}, I_{[\lam]}]\subset I_{[\lam]}.
\end{equation}

Next, we will show that $I_{[\lam]}$ satisfies $I_{[\lam]} I_{[\lam]}\subset I_{[\lam]}.$ We have
\begin{eqnarray}\label{2071}
\nonumber I_{[\lam]} I_{[\lam]}&=&(I_{ 0, [\lam]}\oplus\vv_{[\lam]})
(I_{0,  [\lam]}\oplus\vv_{[\lam]})\\
&\subset&I_{ 0, [\lam]} I_{ 0, \lam]}+I_{ 0, [\lam]} \vv_{[\lam]}+\vv_{[\lam]} I_{ 0, [\lam]}+\vv_{[\lam]} \vv_{[\lam]}.
\end{eqnarray}
It is enough we just have to consider the first summand in (\ref{2071}). For the rest of summands, by a similar way as above, one  can shows
\begin{equation}\label{2072}
I_{0,  [\lam]} \vv_{[\lam]}+\vv_{[\lam]} I_{0,  [\lam]}+\vv_{[\lam]} \vv_{[\lam]}\subset I_{[\lam]}.
\end{equation}
Now, consider the first summand $I_{0,  [\lam]} I_{0,  [\lam]}$ in (\ref{2071}),  we have
\begin{eqnarray}\label{2073}
\nonumber I_{ 0, [\lam]} I_{0,  [\lam]}&=&\left(\sum_{\substack{\mu, \eta\in[\lam]\\ \mu\star\eta=\{0\}}}([\pp_{\mu}, \pp_{\eta}]+\pp_{\mu} \pp_{\eta})\right) \left(\sum_{\substack{\mu', \eta'\in [\lam]\\ \mu'\star\eta'=\{0\}}}([\pp_{\mu'}, \pp_{\eta'}]+\pp_{\mu'} \pp_{\eta'})\right)\\
&\subset&\sum_{\mu, \eta, \mu', \eta'\in [\lam]}( [\pp_{\mu}, \pp_{\eta}] [\pp_{\mu'}, \pp_{\eta'}]+[\pp_{\mu}, \pp_{\eta}](\pp_{\mu'} \pp_{\eta'})\\
\nonumber&+&( \pp_{\mu} \pp_{\eta})[\pp_{\mu'}, \pp_{\eta'}]+(\pp_{\mu} \pp_{\eta})( \pp_{\mu'} \pp_{\eta'})).
\end{eqnarray}

For the  fourth item in (\ref{2073}), By the associativity of product,  we have
\begin{eqnarray*}
\sum_{\substack{\mu, \eta, \mu', \eta'\in [\lam]\\ \mu\star\eta=\{0\}=\mu'\star\eta'}}(\pp_{\mu} \pp_{\eta})( \pp_{\mu'} \pp_{\eta'})&=&\sum_{\substack{\mu, \eta, \mu', \eta'\in [\lam]\\ \mu\star\eta=\{0\}=\mu'\star\eta'}}\left(\pp_{\mu}\left(\pp_{\eta} (\pp_{\mu'} \pp_{\eta'})\right)\right)\\
&\subset&\sum_{\substack{\mu, \eta\in[\lam]\\ \mu\star\eta=\{0\}}}\left(\pp_{\mu}(\pp_{\eta} \pp_0)\right).
\end{eqnarray*}
If $0\neq \pp_{\eta} \pp_0\subset\pp_\tau,$ with $\eta\star 0=\{\tau\},~\tau\in\Lambda_\ss\setminus\{0\}$ then $\{0, \tau\}$ is a connection from $\tau$ to $\eta,$ so $\tau\in[\lam].$ Henc $\pp_{\eta} \pp_0\subset\pp_\tau\subset\vv_{[\lam]}.$ Therefore,
\begin{equation}\label{2073.1}
\sum_{\substack{\mu, \eta, \mu', \eta'\in [\lam]\\ \mu\star\eta=\{0\}=\mu'\star\eta'}}(\pp_{\mu} \pp_{\eta})( \pp_{\mu'} \pp_{\eta'})\subset I_{[\lam]}.
\end{equation}

For the second  item in (\ref{2073}), By  the  Leibniz rule, we have
\begin{eqnarray*}
\sum_{\substack{\mu, \eta, \mu', \eta'\in [\lam]\\ \mu\star\eta=\{0\}=\mu'\star\eta'}}[\pp_{\mu}, \pp_{\eta}](\pp_{\mu'} \pp_{\eta'})\\
&\subset&\sum_{\substack{\mu, \eta\in[\lam]\\ \mu\star\eta=\{0\}}}[\pp_{\mu}, \pp_{\eta}] \pp_0\\
&\subset&\sum_{\substack{\mu, \eta\in[\lam]\\ \mu\star\eta=\{0\}}}\left([\pp_{\mu}, \pp_{\eta}\pp_0]+\pp_{\eta}[\pp_{\mu}, \pp_0]\right).
\end{eqnarray*}
A similar way as above with $\pp_{\eta}\pp_0$ and with $[\pp_{\mu}, \pp_0],$ to get $$\sum_{\substack{\mu, \eta, \mu', \eta'\in [\lam]\\ \mu\star\eta=\{0\}=\mu'\star\eta'}}[\pp_{\mu}, \pp_{\eta}](\pp_{\mu'} \pp_{\eta'})\subset\vv_{[\lam]}.$$  Similarly, $\sum_{\substack{\mu, \eta, \mu', \eta'\in [\lam]\\ \mu\star\eta=\{0\}=\mu'\star\eta'}}( \pp_{\mu} \pp_{\eta})[\pp_{\mu'}, \pp_{\eta'}]\subset\vv_{[\lam]}.$ Therefore,

\begin{equation}\label{2073.2}
\sum_{\substack{\mu, \eta, \mu', \eta'\in [\lam]\\ \mu\star\eta=\{0\}=\mu'\star\eta'}}([\pp_{\mu}, \pp_{\eta}](\pp_{\mu'} \pp_{\eta'})+( \pp_{\mu} \pp_{\eta})[\pp_{\mu'}, \pp_{\eta'}])\subset I_{ [\lam]}.
\end{equation}

For the first item  in (\ref{2073}), taking into account $[\pp_{\mu'},  \pp_{\eta'}]\subset\pp_0$  and   the  Leibniz rule, as above one can get
\begin{equation}\label{2073.3}
\sum_{\substack{\mu, \eta, \mu', \eta'\in [\lam]\\ \mu\star\eta=\{0\}=\mu'\star\eta'}}([\pp_{\mu}, \pp_{\eta}][\pp_{\mu'},  \pp_{\eta'}])\subset I_{ [\lam]}.
\end{equation}

From  Eqs. (\ref{2073.1})-(\ref{2073.3}),  we showed that 
\begin{equation}\label{2074}
 I_{0, [\lam]} I_{0,  [\lam]}\subset I_{[\lam]}.
\end{equation}
Therefore, Eqs. (\ref{2072}) and (\ref{2074})  give us 
\begin{equation}\label{2075}
 I_{[\lam]} I_{[\lam]}\subset I_{[\lam]}.
\end{equation}

Finaly, from Eqs. (\ref{2061}) and (\ref{2075}), we conclude that $[I_{[\lam]}, I_{[\lam]}]+I_{[\lam]} I_{[\lam]}\subset I_{[\lam]},$ so $I_{[\lam]}$ is a (graded) subalgebra of $\pp.$\qed\\

\begin{pro}\label{subalg1}   For any
$\lam, \mu\in\Lam_\ss\setminus\{0\},$   if $[\lam]\neq[\mu]$ then $[I_{[\lam]}, I_{[\mu]}]+I_{[\lam]} I_{[\mu]}=0.$
\end{pro}
\noindent {\bf Proof.} We have
\begin{eqnarray}\label{2}
\nonumber[I_{[\lam]}, I_{[\mu]}]&=&[I_{0,  [\lam]}\oplus\vv_{[\lam]},
I_{0,  [\mu]}\oplus\vv_{[\mu]}]\\
&\subset&[I_{ 0, [\lam]}, I_{0,  [\mu]}]+[I_{0,  [\lam]}, \vv_{[\mu]}]+[\vv_{[\lam]}, I_{0,  [\mu]}]+[\vv_{[\lam]}, \vv_{[\mu]}],
\end{eqnarray}
and also
\begin{eqnarray}\label{3}
\nonumber I_{[\lam]} I_{[\mu]}&=&(I_{0,  [\lam]}\oplus\vv_{[\lam]})
(I_{0,  [\mu]}\oplus\vv_{[\mu]})\\
&\subset&I_{ 0, [\lam]} I_{0,  [\mu]}+I_{0,  [\lam]} \vv_{[\mu]}+\vv_{[\lam]} I_{0,  [\mu]}+\vv_{[\lam]} \vv_{[\mu]}.
\end{eqnarray}

Let us consider the last sammands in (\ref{2}) and (\ref{3}). Suppose that there exist $\lam'\in [\lam]$ and $\eta'\in[\mu]$ such that $[\pp_{\lam'}, \pp_{\mu'}]+\pp_{\lam'} \pp_{\mu'}\neq 0.$  As necessarily $\lam'\star\mu'\neq \emptyset,$ then $\lam', \mu'\in\Lam_\ss\setminus\{0\}.$ So $\{\lam', \mu', \tilde{\lam'}\}$ is a connection from $\lam'$ to $\mu'.$  By the transitivity of the connection relation we have $\mu'\in [\lam],$ a contradiction. Hence,  $[\pp_{\lam'}, \pp_{\mu'}]+\pp_{\lam'} \pp_{\mu'}= 0$ and so
\begin{equation}\label{4}
[\vv_{[\lam]}, \vv_{[\mu]}]+\vv_{[\lam]} \vv_{[\mu]}=0.
\end{equation}

Consider now the second summands in (\ref{2}) and (\ref{3}), we have
\begin{eqnarray*}
[I_{ 0, [\lam]}, \vv_{[\mu]}]+I_{0,  [\lam]} \vv_{[\mu]}&=&\left[\sum_{\substack{\lam_1, \lam_2\in[\lam]\\ \lam_1\star\lam_2=\{0\}}}([\pp_{\lam_1}, \pp_{\lam_2}]+\pp_{\lam_1} \pp_{\lam_2}), \bigoplus_{\mu'\in[\mu]}\pp_{\mu'}\right]\\
&+&\left(\sum_{\substack{\lam_1, \lam_2\in[\lam]\\ \lam_1\star\lam_2=\{0\}}}([\pp_{\lam_1}, \pp_{\lam_2}]+\pp_{\lam_1} \pp_{\lam_2}), \bigoplus_{\mu'\in[\mu]}\pp_{\mu'}\right)\\
&\subset&\sum_{\substack{\lam_1, \lam_2\in[\lam],~ \mu'\in[\mu]\\ \lam_1\star\lam_2=\{0\}}}(\left[[\pp_{\lam_1}, \pp_{\lam_2}], \pp_{\mu'}\right]+[\pp_{\lam_1} \pp_{\lam_2}, \pp_{\mu'}]\\
&+&\left([\pp_{\lam_1}, \pp_{\lam_2}]\pp_{\mu'}+(\pp_{\lam_1}\pp_{\lam_2}) \pp_{\mu'}\right).
\end{eqnarray*}

Suppose there exist  $\lam_1, \lam_2\in [\lam]~~\hbox{with}~~\lam_1\star\lam_2=\{0\}$ and $\mu'\in[\mu]$ such that
\begin{equation*}
[[\pp_{\lam_1}, \pp_{\lam_2}], \pp_{\mu'}]+[\pp_{\lam_1} \pp_{\lam_2}, \pp_{\mu'}]+[\pp_{\lam_1}, \pp_{\lam_2}]\pp_{\mu'}+(\pp_{\lam_1}\pp_{\lam_2}) \pp_{\mu'}\neq 0.
\end{equation*}
The following is divided into four situations to discuss:\\ 

$\bullet$~~If  $[[\pp_{\lam_1}, \pp_{\lam_2}], \pp_{\mu'}]\neq 0.$ By  Jacobi  identity,  we get either $[\pp_{\lam_1}, \pp_{\mu'}]\neq 0$ or $[\pp_{\lam_2}, \pp_{\mu'}]\neq 0.$  In both cases, we have a contradiction, thanks to  equation (\ref{4}). Hence,
 $$
\sum_{\substack{\lam_1, \lam_2\in[\lam],~ \mu'\in[\mu]\\ \lam_1\star\lam_2=\{0\}}}[[\pp_{\lam_1},
\pp_{\lam_2}], \pp_{\mu'}]= 0.
$$

$\bullet$~~If  $[\pp_{\lam_1} \pp_{\lam_2},  \pp_{\mu'}]\neq 0.$ By the Leibniz rule, we have
\begin{equation*}
0\neq[\pp_{\lam_1} \pp_{\lam_2},  \pp_{\mu'}]\subset[\pp_{\lam_1} ,  \pp_{\mu'}] \pp_{\lam_2}+\pp_{\lam_1} [\pp_{\lam_2},  \pp_{\mu'}].
\end{equation*}
So either $[\pp_{\lam_1} ,  \pp_{\mu'}]\neq 0$ or $ [\pp_{\lam_2},  \pp_{\mu'}]\neq 0.$  In both cases, we have a contradiction, thanks to  equation (\ref{4}). Hence, 
$$
\sum_{\substack{\lam_1, \lam_2\in[\lam],~ \mu'\in[\mu]\\ \lam_1\star\lam_2=\{0\}}}[\pp_{\lam_1} 
\pp_{\lam_2}, \pp_{\mu'}]= 0.
$$

$\bullet$~~If  $[\pp_{\lam_1}, \pp_{\lam_2}]\pp_{\mu'}\neq 0.$  By the Leibniz rule, we have
\begin{equation*}
0\neq[\pp_{\lam_1}, \pp_{\lam_2}]\pp_{\mu'}\subset[\pp_{\lam_1} ,  \pp_{\lam_2} \pp_{\mu'}] +\pp_{\lam_2} [\pp_{\lam_1},  \pp_{\mu'}].
\end{equation*}
So either $[\pp_{\lam_1} ,  \pp_{\mu'}]\neq 0$ or $ \pp_{\lam_2} \pp_{\mu'}\neq 0.$  In both cases, we have a contradiction, thanks to  equation (\ref{4}). Hence, 
$$
\sum_{\substack{\lam_1, \lam_2\in[\lam],~ \mu'\in[\mu]\\ \lam_1\star\lam_2=\{0\}}}[\pp_{\lam_1},
\pp_{\lam_2}]\pp_{\mu'}=0.
$$

$\bullet$~~Finally, if    $(\pp_{\lam_1}\pp_{\lam_2}) \pp_{\mu'}\neq 0.$ By asociativity, we have
\begin{equation*}
0\neq(\pp_{\lam_1}\pp_{\lam_2}) \pp_{\mu'}=\pp_{\lam_1}(\pp_{\lam_2} \pp_{\mu'}).
\end{equation*}
So $\pp_{\lam_2} \pp_{\mu'}\neq 0,$ which is  a contradiction (see  equation (\ref{4}). Hence, 
$$
\sum_{\substack{\lam_1, \lam_2\in[\lam],~ \mu'\in[\mu]\\ \lam_1\star\lam_2=\{0\}}}(\pp_{\lam_1}
\pp_{\lam_2}) \pp_{\mu'}= 0.
$$
Therefor,
\begin{equation}\label{6}
[I_{ [\lam], 1}, \vv_{[\mu]}]+I_{ [\lam], 1} \vv_{[\mu]}=0.
\end{equation}
In a similar way, we get
\begin{equation}\label{7}
[\vv_{[\lam]}, I_{ [\mu], 1}]+\vv_{[\lam]} I_{ [\mu], 1}=0.
\end{equation}

Finally, we consider the first summands in  (\ref{2}) and (\ref{3}).  We have
\begin{eqnarray*}
&~& [I_{0,  [\lam]}, I_{0,  [\mu]}]=\\
&~&\left[\sum_{\substack{\lam_1, \lam_2\in[\lam]\\ \lam_1\star\lam_2=\{0\}}}([\pp_{\lam_1}, \pp_{\lam_2}]+\pp_{\lam_1} \pp_{\lam_2}), \sum_{\substack{\mu_1, \mu_2\in[\mu]\\ \mu_1\star\mu_2=\{0\}}}([\pp_{\mu_1}, \pp_{\mu_2}]+\pp_{\mu_1} \pp_{\mu_2})\right]\\
&~&\subset \sum_{\substack{\lam_1, \lam_2\in[\lam],~\mu_1, \mu_2\in[\mu]\\ \lam_1\star\lam_2=\{0\}=\mu_1\star \mu_2}}([[\pp_{\lam_1}, \pp_{\lam_2}], [\pp_{\mu_1}, \pp_{\mu_2}]]+[[\pp_{\lam_1}, \pp_{\lam_2}], \pp_{\mu_1} \pp_{\mu_2}]\\
&+&\left[\pp_{\lam_1} \pp_{\lam_2}, [\pp_{\mu_1}, \pp_{\mu_2}]\right]+[\pp_{\lam_1} \pp_{\lam_2}, \pp_{\mu_1} \pp_{\mu_2}]). 
\end{eqnarray*}
Suppose there exist  $\lam_1, \lam_2\in [\lam]~~\hbox{with}~~\lam_1\star\lam_2=0$ and  $\mu_1, \mu_2\in [\mu]~~\hbox{with}~~\mu_1\star\mu_2=0$  such that
\begin{eqnarray*}
0&\neq&\left[[\pp_{\lam_1}, \pp_{\lam_2}], [\pp_{\mu_1}, \pp_{\mu_2}]\right]+[[\pp_{\lam_1}, \pp_{\lam_2}], \pp_{\mu_1} \pp_{\mu_2}]\\
&+&\left[\pp_{\lam_1} \pp_{\lam_2}, [\pp_{\mu_1}, \pp_{\mu_2}]\right]+[\pp_{\lam_1} \pp_{\lam_2}, \pp_{\mu_1} \pp_{\mu_2}].
\end{eqnarray*}

As above, the  Leibniz rule and Jacobi identity  give us 
\begin{equation*}
[I_{0 [\lam]}, \vv_{[\mu]}]+I_{0,  [\lam]} \vv_{[\mu]}+[\vv_{[\lam]}, I_{0,  [\mu]}]+\vv_{[\lam]} I_{0,  [\mu]}\neq 0,
\end{equation*}
a contradiction either with Eq. (\ref{6}) or with Eq. (\ref{7}). From here
\begin{equation}\label{8}
 [I_{0,  [\lam]}, I_{0,  [\mu]}]=0.
\end{equation}
In a similar manner, one can get
\begin{equation}\label{9}
 I_{0,  [\lam]} I_{0,  [\mu]}=0.
\end{equation}
From Eqs. (\ref{6})- (\ref{9}) and  (\ref{4}), we conclude the result.\qed\\

\begin{thm}\label{main1} The following assertions hold
\begin{itemize}
\item[(1)]  For any $\lambda\in\Lambda_\ss\setminus\{0\},$ the graded subalgebra $I_{[\lam]} =I_{0, [\lam]}\oplus\vv_{[\lam]},$  is a graded  ideal of $\pp.$

\item[(2)] If $\pp$ is gr-simple, then there exists a connection from
$\lam$ to $\mu$ for any $\lam, \mu\in\Lambda_\ss\setminus\{0\}.$ 
\end{itemize}
\end{thm}
\noindent {\bf Proof.} (1) First, we are going to check that
$I_{[\lam]}$ satisfies $[I_{[\lam]}, \pp]\subset I_{[\lam]}.$ We have
\begin{equation}\label{10.}
[I_{[\lam]}, \pp]=[I_{0, [\lam]}\oplus\vv_{[\lam]}, \pp_0\oplus(\bigoplus_{\eta\in\Lambda_\ss\setminus\{0\}}\pp_\eta)].
\end{equation}
Let us consider the product $[\vv_{[\lam]}, \pp_0]$ in Eq. (\ref{10.}) and  uppose there exists $\mu\in[\lam]$ such that $[\pp_\mu, \pp_0]\neq 0.$ We then have $\mu\star 0=\{\eta\}$ with $\eta\in\Lambda_\ss\setminus\{0\}.$ So $\{\mu, o\}$ is a connection from $\mu$ to $\eta$ and then $\eta\in[\lam].$ From here $[\pp_\mu, \pp_0]\subset\pp_\eta\subset\vv_{[\lam]}.$ Hence,

\begin{equation}\label{10}
[\vv_{[\lam]}, \pp_0]\subset I_{[\lam]}.
\end{equation}

Consider the product $[\vv_{[\lam]}, \bigoplus_{\eta\in\Lambda_\ss\setminus\{0\}}\pp_\eta]$  in Eq. (\ref{10.}).  By Propositions \ref{subalg} and \ref{subalg1}, we have 
\begin{equation}\label{11}
[\vv_{[\lam]}, \bigoplus_{\eta\in\Lambda_\ss\setminus\{0\}}\pp_\eta]=\left[\vv_{[\lam]},  \left(\bigoplus_{\eta\in[\lambda]}\pp_\eta\right)\oplus\left(\bigoplus_{\eta\notin[\lambda]}\pp_\eta\right)\right]\subset I_{[\lam]}.
\end{equation}

Consider now the product $[I_{0, [\lam]}, \pp_0]$  in Eq. (\ref{10.}).  We have
\begin{equation*}
[I_{0, [\lam]}, \pp_0]=\left[\sum_{\substack{\mu, \eta\in[\lam]\\ \mu\star\eta=\{0\}}}([\pp_{\mu}, \pp_{\eta}]+\pp_{\mu} \pp_{\eta}, \pp_0\right]\subset\sum_{\substack{\mu, \eta\in[\lam]\\ \mu\star\eta=\{0\}}}([[\pp_\mu, \pp_\eta], \pp_0]+[\pp_\mu\pp_\eta, \pp_0]).
\end{equation*}

Suppose there exist $\mu, \eta\in[\lam]$ with $\mu\star\eta=\{0\}$ such that $0\neq[\pp_\mu\pp_\eta, \pp_0]\subset\pp_\tau$ with $\tau\in\Lam_{\ss}.$ In case  $\tau=0,$ we have $0\neq[\pp_\mu\pp_\eta, \pp_0]\subset\pp_0.$  Now,  by the  Leibniz rule we have either $0\neq \pp_\mu[\pp_\eta, \pp_0]\subset\pp_0$ or $0\neq[\pp_\mu, \pp_0]\pp_\eta\subset\pp_0.$
In the first possibility, we have $\eta\star 0=\{\nu\}$ for some $\nu\in\Lambda_\ss\setminus\{0\}.$ The connection $\{\eta, 0\}$ gives us $\eta\sim\nu$  so $\nu\in[\lam]$ and then 
$
0\neq \pp_\mu[\pp_\eta, \pp_0]\subset\pp_\mu\pp_\nu\subset I_{0, [\lam]}.
$
Similarly in the second possibility we conclude that $0\neq[\pp_\mu, \pp_0]\pp_\eta\subset I_{0, [\lam]}.$ Hence, $\sum_{\substack{\mu, \eta\in[\lam]\\ \mu\star\eta=\{0\}}}[\pp_\mu\pp_\eta, \pp_0]\subset I_{0, [\lam]}.$ We obtain the same
result, by using the  Jacobi identity for the first summand $[[\pp_\mu, \pp_\eta], \pp_0]\neq 0.$  So we can summarize this paragraph by asserting
\begin{equation}\label{12}
[I_{0, [\lam]}, \pp_0]\subset I_{[\lam]}.
\end{equation}

Finally, consider the product $[I_{0, [\lam]}, \bigoplus_{\eta\in\Lambda_\ss\setminus\{0\}}\pp_\eta]$  in Eq. (\ref{10.}).  We have
\begin{eqnarray*}
[I_{0, [\lam]}, \bigoplus_{\nu\in\Lambda_\ss\setminus\{0\}}\pp_\nu]&=&\left[\sum_{\substack{\mu, \eta\in[\lam]\\ \mu\star\eta=\{0\}}}([\pp_{\mu}, \pp_{\eta}]+\pp_{\mu} \pp_{\eta}, \bigoplus_{\nu\in\Lambda_\ss\setminus\{0\}}\pp_\nu\right]\\
&\subset&\sum_{\substack{\mu, \eta\in[\lam],~ \nu\in\Lambda_\ss\setminus\{0\}\\ \mu\star\eta=\{0\}}}([[\pp_\mu, \pp_\eta], \pp_\nu]+[\pp_\mu\pp_\eta, \pp_\nu]).
\end{eqnarray*}
Suppose  there exist $\mu, \eta\in[\lam]$ with $\mu\star\eta=\{0\}$ and $\nu\in\Lambda_\ss\setminus\{0\}$ such that $0\neq[\pp_\mu\pp_\eta, \pp_\nu]\subset\pp_\tau,$ as  necessarily $\tau\in\Lambda_\ss\setminus\{0\}.$ So we have $0\neq[\pp_\mu\pp_\eta, \pp_\nu]\subset\pp_\tau.$   By the  Leibniz rule, we have either $0\neq\pp_\mu[\pp_\eta, \pp_\nu]\subset\pp_\tau$ or $0\neq[\pp_\mu, \pp_\nu]\pp_\eta\subset\pp_\tau.$ In the first possibility, there is a $\delta\in\Lam_{\ss}$ such that $0\neq[\pp_\eta, \pp_\nu]\subset\pp_\delta$ with $\eta\star\nu=\{\delta\}.$ So $\{\eta, \delta\}$ is a connection from $\eta$ to $\tau$ and then $\tau\in[\lam].$ Hence $\pp_\mu[\pp_\eta, \pp_\nu]\subset\pp_\mu\pp_\delta\subset\vv_{[\lam]}.$ Similarly in the second possibility we conclude that  $[\pp_\mu, \pp_\nu]\pp_\eta\subset\vv_{[\lam]}.$ Therfore, $\sum_{\substack{\mu, \eta\in[\lam]\\ \nu\in\Lambda_\ss\setminus\{0\}}}[\pp_\mu\pp_\eta, \pp_\nu]\subset\vv_{[\lam]}.$ We obtain the same  result, by using the  Jacobi identity for the first summand $[[\pp_\mu, \pp_\eta], \pp_\nu]\neq 0.$ From here
\begin{equation}\label{13}
[I_{0, [\lam]}, \bigoplus_{\eta\in\Lambda_\ss\setminus\{0\}}\pp_\eta]\subset I_{[\lam]}.
\end{equation}

From Eqs. (\ref{10})-(\ref{11}), we get 
\begin{equation*}
[I_{[\lam]}, \pp]\subset I_{[\lam]}.
\end{equation*}
Next, we will check that $I_{[\lam]} \pp\subset I_{[\lam]}.$ We have
\begin{equation}\label{14}
I_{[\lam]}\pp=(I_{0, [\lam]}\oplus\vv_{[\lam]})\left(\pp_0\oplus(\bigoplus_{\eta\in[\lam]}\pp_\eta)\oplus(\bigoplus_{\eta\notin[\lam]}\pp_\eta)\right).
\end{equation}
Let us consider the product  $I_{0, [\lam]}\pp_0$ in Eq. (\ref{14}).  By using the asociatitvity and the Leibniz rule as the proof of Proposition \ref{subalg} we have 
\begin{equation*}
I_{0, [\lam]}\pp_0=\sum_{\substack{\mu, \eta\in[\lam]\\ \mu\star\eta=\{0\}}}([\pp_{\mu}, \pp_{\eta}]+\pp_{\mu} \pp_{\eta})\pp_0\subset\sum_{\substack{\mu, \eta\in[\lam]\\ \mu\star\eta=\{0\}}}\left([\pp_{\mu}, \pp_{\eta}]\pp_0+(\pp_{\mu} \pp_{\eta})\pp_0\right)\subset I_{0, [\lam]}.
\end{equation*}
We also have $\vv_{[\lam]}\pp_0\subset\vv_{[\lam]}.$ Thus $(I_{0, [\lam]}\oplus\vv_{[\lam]})\pp_0\subset I_{[\lam]}.$ From here and  Propositions \ref{subalg} and \ref{subalg1}, we can assert 
$$
I_{[\lam]}\pp=(I_{0, [\lam]}\oplus\vv_{[\lam]})\left(\pp_0\oplus(\bigoplus_{\eta\in[\lam]}\pp_\eta)\oplus(\bigoplus_{\eta\notin[\lam]}\pp_\eta)\right)\subset I_{[\lam]}.
$$
In a similar way we get $\pp I_{[\lam]}\subset I_{[\lam]}$  and so  $ I_{[\lam]}$ is an ideal of $\pp.$\\

(2) The simplicity of $\pp$ implies that $\pp=I_{[\lam]}$  for some $\lam\in\Lambda_\ss\setminus\{0\}.$ Hence, $[\lam]=\Lambda_\ss\setminus\{0\}$ and so any pair of elements in $\Lambda_\ss\setminus\{0\}$ is connected. \qed\\

\begin{thm}\label{main2} A set-graded non-commutative Poisson algebra $\pp$ decompose as
$$
\pp=\uu\oplus\sum_{[\lam]\in(\Lambda_\ss\setminus\{0\})/\sim}I_{[\lam]},
$$
where $\uu$ is a linear space complement of $\span_{\bbbf}\{ [\pp_{\mu}, \pp_{\eta}]+\pp_{\mu}\pp_{\eta}  : \mu, \eta\in[\lam]\}$ in $\pp_0$ and
 any $I_{[\lam]}$ is one of the graded ideals of $\pp$ described in Theorem \ref{main1}-(1), satisfying $[I_{[\lam]}, I_{[\mu]}]+I_{[\lam]} I_{[\mu]}=0,$ whenever $[\lam]\neq[\mu].$
\end{thm}
\noindent {\bf Proof.} We have $I_{[\lam]}$  well-defined and, by Theorem  \ref{main1}-(1), a graded ideal of $\pp.$  Now, by
considering a linear complement $\u$ of $\span_{\bbbf}\{ [\pp_{\mu}, \pp_{\eta}]+\pp_{\mu}\pp_{\eta}  : \mu, \eta\in[\lam]\}$  in $\pp_0,$  we have
$$
\pp=\pp_0\oplus(\bigoplus_{\lam\in\Lambda_\ss\setminus\{0\}}\pp_\lam)=\uu\oplus\sum_{[\lam]\in(\Lambda_\ss\setminus\{0\})/\sim}I_{[\lam]}.
$$
Finally, Proposition \ref{subalg1} gives us $[I_{[\lam]}, I_{[\mu]}]+I_{[\lam]} I_{[\mu]}=0,$ whenever $[\lam]\neq[\mu].$\qed\\

In case it is not distinguished any element $0$ in the support of
the grading, that is $0=\emptyset,$ we have as an immediate consequence
of Theorem  \ref{main1} the following result;
\begin{cor}\label{10001}
If $0=\emptyset,$ then
$$
 \pp=\bigoplus_{[\lam]\in(\Lambda_\ss\setminus\{0\})/\sim}I_{[\lam]},
$$
where  any $I_{[\lam]}$ is one of the graded ideals of $\pp$ described in Theorem \ref{main1}-(1), satisfying $[I_{[\lam]}, I_{[\mu]}]+I_{[\lam]} I_{[\mu]}=0,$ whenever $[\lam]\neq[\mu].$\qed
\end{cor}

Let us denote by $Z(\pp)$ the centre  of $\pp,$ that is,
$Z(\pp)=\{x\in\pp~~:~~[x, \pp]+x \pp+\pp x=0\}.$\\

\begin{DEF} Let $\pp$ be a set-graded non-commutative Poison algebra. We say that $\pp_0$  is  tight whence 
$$
\pp_0=\{0\}~~~\hbox{or}~~~\pp_0=\sum_{\substack{\mu, \eta\in\Lambda_\ss\setminus\{0\}\\ \lam\star\mu=\{0\} }}([\pp_\lam, \p_\mu]+\pp_\lam, \p_\mu).
$$
\end{DEF}

\begin{cor}\label{2.17}  If $Z(\pp)=0$ and  $\pp_0$  is  tight.  Then $\pp$ is the direct sum of the graded ideals
given in Theorem \ref{main1}-(1),
$$
\pp=\bigoplus_{[\lam]\in(\Lambda_\ss\setminus\{0\})/\sim}I_{[\lam]},
$$
with $[I_{[\lam]}, I_{[\mu]}]+I_{[\lam]} I_{[\mu]}=0,$ whenever $[\lam]\neq[\mu].$
\end{cor}
\noindent {\bf Proof.} .Since $\pp_0$  is  tight,  it is clear that $$\pp=\sum_{[\lam]\in(\Lam_\ss\setminus\{0\})/\sim}I_{[\lam]}.$$

 For the direct character, take some $$x\in I_{[\lam]}\cap\sum_{\substack{[\mu]\in(\Lambda_\ss\setminus\{0\})/\sim\\  [\lam]\neq[\mu]}}I_{[\mu]}.$$ From $x\in I_{[\lam]}$ and the fact $[I_{[\lam]}, I_{[\mu]}]+I_{[\lam]} I_{[\mu]}+ I_{[\mu]}I_{[\lam]}=0,$ if $[\lam]\neq[\mu],$  we get
\begin{equation}\label{17}
\left[x, \sum_{\substack{[\mu]\in(\Lambda_\ss\setminus\{0\})/\sim\\  [\lam]\neq[\mu]}}I_{[\mu]}\right]+x\left (\sum_{\substack{[\mu]\in(\Lambda_\ss\setminus\{0\})/\sim\\  [\lam]\neq[\mu]}}I_{[\mu]}\right)+\left(\sum_{\substack{[\mu]\in(\Lambda_\ss\setminus\{0\})/\sim\\  [\lam]\neq[\mu]}}I_{[\mu]}\right) x=0.
\end{equation}
In the other hand, since $x\in\sum_{\substack{[\mu]\in(\Lambda_\ss\setminus\{0\})/\sim\\  [\lam]\neq[\mu]}}I_{[\mu]}$ and the same above fact implies that 
\begin{equation}\label{18}
[x, I_{[\lam]}]+x I_{[\lam]}+I_{[\lam]} x=0.
\end{equation}
Now, Eqs. (\ref{17}) and (\ref{18}) give us  $x\in Z(\pp)$ and so $x=0.$
Hence, $\pp=\bigoplus_{[\lam]\in(\Lambda_\ss\setminus\{0\})/\sim}I_{[\lam]},$   as desired\qed\\

\section{The graded simple components} \setcounter{equation}{0}\

In this section, we will study the  simplicity  of  set-graded non-commutative  Poisson   algebras and interested in studing   under which conditions a set-graded non-commutative  Poisson   algebra $\pp$ decomposes as the direct sum of the family of its
gr-simple ideals. We begin by introducing the concepts of maximal length and
$\Lam_{\ss}$-multiplicativity in the setup of set-graded non-commutative  Poisson   algebras  with a set grading in a similar way as in the frameworks of set-graded Lie algebras, set-graded Lie superalgebras, split non-commutative Poisson algebras and group-graed Poisson color algebras etc. (see \cite{C0,  C1, Kh1, Kh2} for discussions and examples on these concepts). From now on,  for any $\lam\in\Lambda_{\ss}$ we will denote $\pp_{\tilde{\lam}}=\{0\}.$\\

\begin{DEF} We say that a set-graded non-commutative  Poisson   algebra $\pp$  is $\Lam_{\ss}$-multiplicative if given $\lam, \mu\in\Lam_{\ss}$ such that $\lam\in\mu\star s$ for some $s\in\Lambda_\ss\dot{\cup}\widetilde{\Lambda_\ss}$ then
$$
\pp_\lam\subset[\pp_{\mu}, \pp_s+\pp_{\tilde{s}}]+\pp_{\mu}(\pp_s+\pp_{\tilde{s}}).
$$
\end{DEF}

\begin{DEF}  A set-graded non-commutative  Poisson   algebra $\pp$ is of maximal length if for any $\lam\in\Lambda_\ss\setminus\{0\}$ we have $\dim \pp_{\lam}=1.$
\end{DEF}

\begin{lem}\label{4.1} Let $\pp$ be a centerless  set-graded non-commutative Poisson  algebra of maximal length and  with $\pp_0$ tight.  If $I$ is an  ideal of $\pp$ such that $I\subset\pp_0$ then $I=\{0\} .$
\end{lem}
\noindent {\bf Proof.} Suppose there exists a nonzero graded ideal $I=\bigoplus_{\lam\in\Lam_\ss}I_\lam$ of $\pp$  such that $I\subset\pp_0.$ The fact that $I_\lam=I\cap\pp_\lam,$ for any $\lam\in\Lam_\ss$ and the maximal length of $\pp$  get us 
\begin{equation}\label{41}
I=(I\cap\pp_0)\oplus\left(\bigoplus_{\lam\in\Lambda_\ss\setminus\{0\}}(I\cap\pp_\lam)\right).
\end{equation}
Now, given any $\lam\in\Lambda_\ss\setminus\{0\},$ taking into account $I\subset\pp_0, $ we have
$$
[I, \pp_\lam]+I \pp_\lam+\pp_\lam I\subset \pp_0\cap\pp_\mu,
$$
with $0\star\lam=\{\mu\}$ for some $\mu\in\Lambda_\ss\setminus\{0\}.$ Hence,
\begin{equation}\label{42}
[I, \pp_\lam]+I \pp_\lam+\pp_\lam I=0,~~\forall\lam\in\Lambda_\ss\setminus\{0\}.
\end{equation}
Since $Z(\pp)=0$ and $\pp_0$ is tight, for any $0\neq x\in I$ we have there exist $\lam, \mu\in\Lambda_\ss\setminus\{0\}$ such that
$$
[x, [\pp_\lam, \pp_\mu]]+x (\pp_\lam \pp_\mu)+(\pp_\lam \pp_\mu)x\neq 0.
$$
By the Jacobi identity and associativity, we have either $[x, \pp_\lam]+x \pp_\lam+\pp_\lam x\neq 0$ or $[x, \pp_\mu]+x \pp_\mu+\pp_\mu x\neq 0$ which is a contradiction with Eq. (\ref{42}). Therefore, we conclude $I=\{0\}.$\qed\\

\begin{thm}\label{main4} Let $\pp$ be a a centerless $\Lam_{\ss}$-multiplicative set-graded non-commutative Poisson  algebra of maximal length and  with $\pp_0$ tight. Then $\pp$ is gr-simple if and only if it has all of the non-zero elements in $\Lam_\ss$ connected. 
\end{thm}
\noindent {\bf Proof.} The first implication is Theorem \ref{main1}-(2).  To prove the converse, consider $I=\bigoplus_{\lam\in\Lam_\ss}I_\lam$ a non-zero ideal of $\pp.$ Taking into account Eq.(\ref{41}),  we can write
\begin{equation}\label{44}
 I=(I\cap\pp_0)\oplus(\bigoplus_{\lam\in\Lam^I_{\ss}} \pp_{\lam}),
\end{equation}
where  $\Lam^I_{\ss} :=\{\lam\in\Lambda_\ss\setminus\{0\}~:~I\cap\pp_\lam\neq \{0\}\},$ also being $\Lam^I_{\ss} \neq\emptyset,$  as a consequence of Lemma \ref{4.1}.  Hence, we may choose $\lam_0\in\Lam_\ss^I$  being so
\begin{equation}\label{45}
0\neq \pp_{\lam_0}\subset I.
\end{equation}
Now, let us take any $\mu\in\Lambda_\ss\setminus\{0\}.$  The fact that $\lam_0$ is connected to $\mu,$
gives us a connection $\{\lam_1, \lam_2, \lam_3, ..., \lam_k\}\subset\Sigma_{\Lam}\subset\Lambda_\ss\dot{\cup}\widetilde{\Lambda_\ss},$
satisfying the following conditions;\\

If $k=1:$
\begin{itemize}
\item[(1)] $\lam_1=\lambda_0=\mu.$
\end{itemize}

If $k\geq2:$
\begin{itemize}
\item[(1)] $\lam_1\in\{\lambda, \tilde{\lambda}\}$,
\item[(2)] $\psi(\{\lam_1\}, \lam_2)\neq\emptyset,\\
\psi(\psi(\{\lam_1\}, \lam_2), \lam_3)\neq\emptyset,\\
\psi(\psi(\psi(\{\lam_1\}, \lam_2), \lam_3), \lam_4)\neq\emptyset,\\
...\\
\psi(\psi(...(\psi(\{\lam_1\}, \lam_2), ...), \lam_{k-2}),
\lam_{k-1})\neq\emptyset$.\\
\item[(3)] $\mu\in\psi(\psi(...(\psi(\{\lam_1\}, \lam_2), ...), \lam_{k-2}),
\lam_{k-1}), \lam_k)$.
\end{itemize}
 Consider $\psi(\{\lam_1\}, \lam_2)\neq\emptyset$
and so $\psi(\{\lam_1\}, \lam_2)\cap\Lambda_\ss\neq\emptyset$ (see Remark
\ref{1111}).  Hence, for any $\mu_1\in\psi(\{\lam_1\},
\lam_2)\cap\Lambda_\ss,$ we have $\mu_1\in\lam_1\star\lam_2.$ The fact $\lam_1\in\{\lambda, \tilde{\lambda}\}$ gives us either $\mu_1\in\lam_0\star\lam_2$ or $\mu_1\in\tilde{\lam_0}\star\lam_2=\lam_2\star\tilde{\lam_0}$ with necessarily $\lam_2\in\Lam_\ss$ in the second possibility. Now, by $\Lam_\ss-$multiplicativity and    maximal length of $\pp$ as  consequence of
Eq. (\ref{45}), we get either
\begin{equation*}
0\neq\pp_{\mu_1}=[\pp_{\lam_0}, \pp_{\lam_2}+\pp_{\tilde{\lam_2}}]+\pp_{\lam_0} (\pp_{\lam_2}+\pp_{\tilde{\lam_2}})\subset I.
\end{equation*}
or $0\neq\pp_{\mu_1}=[\pp_{\lam_0}, \pp_{\lam_2}]+\pp_{\lam_0}\pp_{\lam_2}\subset I.$
From here  we assert that
\begin{equation}\label{47}
\bigoplus_{\lam\in\psi(\{\lam_1\},
\lam_2)\cap\Lambda_\ss}\pp_\lam\subset I.
\end{equation}
Next, again Remark \ref{1111} shows that $\psi(\psi(\{\lam_1\}, \lam_2),
\lam_3) \cap\Lambda_\ss\neq\emptyset.$  Given any
$\mu_2\in\psi(\psi(\{\lam_1\}, \lam_2), \lam_3) \cap\Lambda_\ss,$ we have $\mu_2\in\eta\star\lam_3$ for some $\eta\in\psi(\{\lam_1\}, \lam_2).$  Now, by $\Lam_\ss-$multiplicativity and    maximal length of $\pp$ as  consequence of
Eq. (\ref{47}), we get either
\begin{equation*}
0\neq\pp_{\mu_2}=[\pp_{\eta}, \pp_{\lam_3}+\pp_{\tilde{\lam_3}}]+\pp_{\eta} (\pp_{\lam_3}+\pp_{\tilde{\lam_3}})\subset I.
\end{equation*}
or $\eta\in\psi(\{\lam_1\}, \lam_2) \cap\Lambda_\ss$ and $\lam_3\in\Lam_\ss,$ taking into account $\mu_2\in\lam_3\star\eta$ and $\tilde{\mu_2}\in\psi(\{\lam_1\}, \lam_2) \cap\Lambda_\ss$ we get
\begin{equation*}
0\neq\pp_{\mu_2}=[\pp_{\lam_3}, \pp_{\tilde{\eta}}]+\pp_{\lam_3}\pp_{\tilde{\eta}}\subset I.
\end{equation*}
From here  we assert that
\begin{equation*}
\bigoplus_{\lam\in\psi(\psi(\{\lam_1\}, \lam_2),
\lam_3) \cap\Lambda_\ss}\pp_\lam\subset I.
\end{equation*}
Following this process with the connection $\{\lam_1, \lam_2, \lam_3, ..., \lam_k\}\subset\Sigma_{\Lam}\subset\Lambda_\ss\dot{\cup}\widetilde{\Lambda_\ss},$  we obtain that
\begin{equation*}
\bigoplus_{\lam\in\psi(\psi(...(\psi(\{\lam_1\}, \lam_2), ...), \lam_{k-2}),
\lam_{k-1}), \lam_k)\cap\Lambda_\ss}\pp_\lam\subset I.
\end{equation*}
Taking now into account $\mu\in\psi(\psi(...(\psi(\{\lam_1\}, \lam_2), ...), \lam_{k-2}),
\lam_{k-1}), \lam_k)\cap\Lambda_\ss$ we conclude that $\pp_\mu\subset I$ and so
\begin{equation}\label{48}
\bigoplus_{\mu\in\Lambda_\ss\setminus\{0\}}\pp_\mu\subset I.
\end{equation}
Finally, the fact that $\pp_0$ is tight together with Eq. (\ref{48}) gives us $\pp_0\subset I$ and so $I=\pp.$  That
is, $\pp$ is gr-simple.\qed\\

\begin{thm}\label{main4} Let $\pp$ be a a centerless $\Lam_{\ss}$-multiplicative set-graded non-commutative Poisson  algebra of maximal length and  with $\pp_0$ tight. Then
$$
\pp=\bigoplus_{[\lam]\in(\Lambda_\ss\setminus\{0\})/\sim}I_{[\lam]},
$$
where any $I_{[\lam]}$ is a simple ideal having  all of its lements different to $0$  in its support connected.
\end{thm}
\noindent {\bf Proof.} By corollary \ref{2.17}, $\pp=\bigoplus_{[\lam]\in(\Lambda_\ss\setminus\{0\})/\sim}I_{[\lam]},$ is the direct sum of the ideals
\begin{eqnarray*}
I_{[\lam]} &=&I_{0, [\lam]}\oplus\vv_{[\lam]}\\
&=&\left(\span_{\bbbf}\{ [\pp_{\mu}, \pp_{\eta}]+\pp_{\mu}\pp_{\eta}  : \mu, \eta\in[\lam]\}\cap\pp_0\right)\oplus(\bigoplus_{\mu\in[\lam]}\pp_\mu),
\end{eqnarray*}
 having any $I_{[\lam]}$  its
support, $\Lam_{I_{[\lam]}}=[\lam]$ (connected through the elements contained in $[\lam]\dot{\cup}[\lam]$). In order to apply Theorem \ref{main4} to each $I_{[\lam]},$
observe that the fact $\Lam_{I_{[\lam]}}=[\lam],$ gives us easily that $\Lam_{I_{[\lam]}}$  has all of its lements different to $0$  in its support connected. We also
have that any of the $I_{[\lam]}$ is $\Lam_{I_{[\lam]}}$-multiplicative as consequence of the $\Lam_\ss$-tmultiplicativity of $\pp.$ Clearly $I_{[\lam]}$  is of maximal length. The $0$-homogeneous component of $I_{[\lam]}$
is equal to $I_{0, [\lam]}$  and so tight by construction. We also have  $Z_{I_{[\lam]}}(I_{[\lam]})=0,$ the centre of $I_{[\lam]}$ in itself, as consequence of  $[I_{[\lam]}, I_{[\mu]}]+I_{[\lam]} I_{[\mu]}=0$ if $[\lam]\neq[\mu]$ (see Proposition \ref{subalg1}) and $Z(\pp)=0.$ 
We can apply Theorem \ref{main4} to any $I_{[\lam]}$ so as to conclude that
$I_{[\lam]}$ is gr-simple. It is clear that the decomposition $\pp=\bigoplus_{[\lam]\in(\Lambda_\ss\setminus\{0\})/\sim}I_{[\lam]},$ satisfies the assertions of the theorem.\qed

\end{document}